\documentclass[11pt]{article}

\usepackage[margin=1.3in]{geometry}

\usepackage{amsmath}
\usepackage{amssymb}
\usepackage{amsthm}

\newtheorem{theorem}{Theorem}
\newtheorem{lemma}{Lemma}

\newtheorem*{theorem*}{Theorem}

\theoremstyle{definition}
\newtheorem*{definition}{Definition}
\newtheorem*{corollary}{Corollary}
\newtheorem*{question}{Question}

\newcommand{\Z}{\mathbb{Z}}

\title{Structural Results on the Unfair 0-1 Polynomial Conjecture}
\author{Emina Hasanbegovi\'c}

\begin{document}

\maketitle

\begin{abstract}
A $0$--$1$ polynomial is a polynomial whose coefficients are all
either $0$ or $1$. 
The unfair \(0\)--\(1\) polynomial conjecture states that if \(f\) is
a $0$--$1$ polynomial and
$
    f(x)=g(x)h(x)
$
is a factorisation into non-constant polynomials with non-negative real
coefficients, then \(g\) and \(h\) must themselves be \(0\)--\(1\)
polynomials. We call \(f\) \emph{fair} if every such factorisation has
this property, and \emph{unfair} otherwise. In this paper, we study
structural properties of possible unfair factorisations and obtain
conditions which guarantee fairness. In particular, we prove that if
\[
    f(x)=x^{\,n-\deg w}w(x)+v(x),
\]
where \(v\) and \(w\) are \(0\)--\(1\) polynomials and \(v\) is fair,
then \(f\) is fair whenever
\[
    n\geq\max\{\deg w+2\deg v,\,2\deg w\}.
\]
\end{abstract}

\section{Introduction}
The unfair \(0\)--\(1\) polynomial
conjecture asks whether every factorisation
$
    f(x)=g(x)h(x)
$
of a \(0\)--\(1\) polynomial \(f\), that is a polynomial whose coefficients all belong to \(\{0,1\}\), into polynomials \(g\) and \(h\)
with non-negative real coefficients necessarily has both \(g\) and
\(h\) themselves \(0\)--\(1\) polynomials.
We call such a factorisation \emph{fair} if both \(g\) and \(h\) are
\(0\)--\(1\) polynomials, and \emph{unfair} otherwise. Similarly, we call \(f\)
\emph{fair} if every factorisation of \(f\) into polynomials with
non-negative real coefficients is fair, and \emph{unfair} otherwise.

This terminology reflects an equivalent probabilistic formulation of
the conjecture. Let \(X\) and \(Y\) be independent discrete random
variables with finite support. If \(X+Y\) is uniformly distributed on
its support, must \(X\) and \(Y\) themselves be uniformly distributed
on their respective supports? Equivalently, if two dice have a fair
sum, must the individual dice themselves be fair?

To the best of the author's knowledge, the question in this general
form was first posed by G.~Letac around 1969~\cite{Ghidelli}. An important
special case, however, predates the general conjecture. In 1937,
Krasner and Ranulac~\cite{KrasnerRanulac} proved the conjecture for
polynomials of the form
\[
    1+x+\cdots+x^{n-1}+x^n
\]
for arbitrary \(n\). Their argument is elementary and extends more
generally to palindromic \(0\)--\(1\) polynomials; see also
Ghidelli~\cite{Ghidelli}. 

More recent work has focused on restricting the possible forms of
unfair factors. Ghidelli~\cite{Ghidelli} studied the first nontrivial
candidate family
\[
    1+ax^2+x^5, \qquad 0<a<1,
\]
and proved that no polynomial of this form can occur as a factor in
an unfair factorisation. His proof combines linear recurrences,
computer algebra, and analytic and numerical methods.
Hare~\cite{Hare} subsequently developed computational methods based
on coefficient constraints, recursive case analysis, and symbolic
techniques to study possible unfair factors of higher degree. Using these methods, he ruled out most candidate factors of degree up to \(15\), eliminating all but \(975\)
of \(7\,141\,686\) possible cases.

Further computational evidence also supports the conjecture. As
reported by Ghidelli~\cite{Ghidelli}, the conjecture was verified for
all \(0\)--\(1\) polynomials of degree at most \(12\) by Sil,
subsequently up to degree \(26\) by Max Alekseyev, and up to degree
\(32\) by Peter Mueller. 

The aim of this paper is to establish structural restrictions on
possible unfair factorisations and thereby narrow the form that a
counterexample to the conjecture could take. In particular, we study
the coefficient structure of a factorisation \(f=gh\), obtain a
condition involving palindromic supports, and, motivated by work on
the irreducibility of polynomials with large
gaps~\cite{SawinShustermanStoll,FilasetaMatthews}, consider
\(0\)--\(1\) polynomials whose non-zero coefficients are concentrated
in two separated blocks. For the latter, we obtain the following
result.

\begin{theorem*}
Let \(v(x)\) and \(w(x)\) be non-zero \(0\)--\(1\) polynomials and suppose $v(x)$ is fair. If
\[
    n\geq \max\{\deg w+2\deg v,\,2\deg w\}
\]
then the polynomial
\[
    f(x)=x^{\,n-\deg w}w(x)+v(x)
\]
is fair.
\end{theorem*}

% ============================================================
\section{Structure of Factorisations}
% ============================================================

The coefficient-wise arguments in this section are closely related to techniques used by Ghidelli [3] and Hare [1], in particular Hare’s trinary framework. In this paper, we use a related viewpoint to study the coefficients of both factors simultaneously through the pairs $(b_i,c_i)$.

We begin by establishing some basic properties of such factorisations. Suppose
\[
    f(x)=g(x)h(x),
\]
and 
\[
g(x)=b_rx^r+\cdots+b_1x+b_0
\]
\[
h(x)=c_sx^s+\cdots+c_1x+c_0
\]
where \(g,h\) are polynomials with non-negative coefficients and $n=r+s$. From simple observations, we can conclude that the coefficients $b_i$ and $c_i$ must lie in the interval $[0,1]$. If we assume that the factorisation is unfair then it is easily seen that $b_0,c_0,b_r$ and $c_s$ must be equal to $1$. Additionally, if $g$ is a $0$--$1$ polynomial then $h$ must be as well.

If we rewrite $f$ in the following way
\[
    f(x)
    =x^{n_k}+x^{n_{k-1}}+\cdots+x^{n_1}+1
\]
that is,
\[
    a_{n_i}=1
    \qquad (i=0,1,\ldots,k),
\]
with the convention \(n_0=0\), and
$
    a_j=0
    \text{ for } j\notin\{n_0,n_1,\ldots,n_k\},
$
we get the following lemma.

\begin{lemma}
Only the coefficients \(b_{n_i}\) and \(c_{n_i}\), for
\(i=0,1,\ldots,k\), can be non-zero.
\end{lemma}

\begin{proof}
Suppose that \(b_m\neq 0\) for some
$
    m\notin\{n_0,n_1,\ldots,n_k\}.
$
Since \(b_0=c_0=1\), we have
\[
    b_mc_0=b_m>0.
\]
Hence the coefficient \(a_m\) must be positive. Since \(f\) is a
\(0\)--\(1\) polynomial, this implies \(a_m=1\). Therefore
\(m=n_i\) for some \(i=0,1,\ldots,k\), a contradiction. The argument for \(c_m\neq0\) is identical.
\end{proof}

\begin{lemma}\label{lm:no11}
For a fixed \(i>0\), if \(b_i=1\), then \(c_i=0\), and vice versa.
\end{lemma}

\begin{proof}
Suppose \(b_i=1\). Since \(b_0=c_0=1\),
\[
    a_i
    \ge b_ic_0+b_0c_i
    =b_i+c_i
    =1+c_i.
\]
Since \(a_i\in\{0,1\}\), the only possibility is
$
    a_i=1, 
    c_i=0.
$
The other case is symmetric.
\end{proof}

This leaves us with several possibilities for the pairs
\((b_i,c_i)\). Following the notation introduced by Kevin G. Hare [1], we write \(*\) for a coefficient strictly between \(0\) and 1, that is
\[
    0<*<1.
\]
Thus, for \(i>0\), the pair \((b_i,c_i)\) can be one of the following:

\begin{itemize}
    \item A \textbf{zero} pair: $(b_i,c_i)=(0,0)$.
    
    \item A \textbf{whole} pair: $(b_i,c_i)=(1,0)$ or $(0,1)$.
    
    \item A \textbf{half} pair: $(b_i,c_i)=(*,0)$ or $(0,*)$.
    
    \item A \textbf{double} pair: $(b_i,c_i)=(*,*)$.
\end{itemize}

For \(i=0\), we have \(b_0=c_0=1\), and by the preceding lemma this
is the only index at which both coefficients can equal \(1\).

\begin{lemma}
Let \((b_i,c_i)\) be the first half pair. Then there is some $j<i$ such that \((b_j,c_j)\) is a double pair.
\end{lemma}

\begin{proof}
Suppose there are no double pairs before \((b_i,c_i)\), and assume
without loss of generality that 
$
    (b_i,c_i)=(*,0).
$
Then
\[
\begin{aligned}
    a_i
    &=\sum_{j+k=i} b_jc_k 
    =b_i+c_i
      +\sum_{\substack{j+k=i\\j<i,\;k<i}}b_jc_k \\
    &=*+z,
\end{aligned}
\]
where \(z\in\Z\), because all earlier coefficients \(b_j,c_k\) are
either \(0\) or \(1\). Since \(a_i\in\Z\), it follows that \(*\) must also be an integer,
contradicting \(0<*<1\).
\end{proof}

\begin{theorem}\label{thm:doubletrouble}
If there are no double pairs \((b_{n_i},c_{n_i})\) for \(i>0\),
then the factorisation must be fair.
\end{theorem}

\begin{proof}
Suppose there are no double pairs. By the previous lemma, there are no
half pairs either. Hence all coefficients of \(g\) and \(h\) are either
\(0\) or \(1\), and so the factorisation is fair.
\end{proof}

This result motivates us to look only at the set $(n_1,...,n_k)$ and consider for which values of $n_i$ we can have $(b_{n_i},c_{n_i})$ be a double pair. 

\begin{corollary}
If $2n_i \notin \{n_1,\ldots,n_k\}$
for every 
$ i=1,\ldots,k,$
then \(f\) is fair.
\end{corollary}

\begin{proof}
For every \(i>0\), the assumption gives \(a_{2n_i}=0\). Hence
\[
    b_{n_i}c_{n_i}=0,
\]
so \((b_{n_i},c_{n_i})\) cannot be a double pair. Therefore there are
no double pairs, and Theorem~\ref{thm:doubletrouble} implies that \(f\) is fair.
\end{proof}

% ============================================================
\section{Symmetry}
% ============================================================

\begin{definition}
Let \(f\) be a polynomial of degree \(n\). We define its
\textit{reciprocal polynomial} \(f^*\) by
\[
    f^*(x):=x^n f(x^{-1}).
\]
\end{definition}
Note that
$
    (f^*)^*=f.
$

\begin{lemma}
Let \(f\) be a \(0\)--\(1\) polynomial. Then \(f\) is fair if and
only if \(f^*\) is fair.
\end{lemma}

\begin{proof}
Suppose \(f^*\) is fair. Let \(f=gh\), where \(g\) and \(h\) have
non-negative coefficients. Taking reciprocals, we have \(f^*=g^*h^*\).
Since \(f^*\) is fair, \(g^*\) and \(h^*\) must be \(0\)--\(1\)
polynomials, so \(g\) and \(h\) are themselves \(0\)--\(1\) polynomials.
Therefore \(f\) is fair.

Conversely, if \(f\) is fair, applying the same argument to \(f^*\),
together with \((f^*)^*=f\), shows that \(f^*\) is fair.
\end{proof}

Considering that this is a symmetric condition, it is interesting to ask what happens in the case $f=f^*$. Remarkably, the conjecture was proved in the palindromic case by Krasner and
Ranulac~\cite{KrasnerRanulac} as mentioned in the introduction. The following two results can also be shown.

\begin{theorem} The polynomial $g(x)$ is a factor of an unfair polynomial if and only if $g^*(x)$ is a factor of an unfair polynomial.
\end{theorem}
\begin{proof}
See Kevin G. Hare~\cite{Hare}.
\end{proof}

For the next result, we call a set
\(\{m_0,m_1,\ldots,m_q\}\), where
\(0=m_0< m_1<\cdots<m_q\), \textit{palindromic} if
\[
    m_i+m_{q-i}=m_q
    \qquad\text{for all } i=0,1,\ldots,q.
\]

\begin{theorem}
If the subset of $(n_0,n_1,...,n_k)$ for which $b_{n_i}\neq 0$ forms a palindromic set, then $g(x)$ must be a \(0\)--\(1\) polynomial.
\end{theorem}
\begin{proof}
Let the subset of $(n_0,n_1,...,n_k)$ for which $b_{n_i}\neq 0$ be $(m_0,m_1,...,m_q)$. Then
\[
g(x)=1+b_{m_1}x^{m_1}+...+b_{m_{q-1}}x^{m_{q-1}}+x^{m_q}
\]
where $m_i+m_{q-i}=m_q$ for each $i=0,1,...,q$. Consider $a_{m_q}$:  
\[
a_{m_q}=b_{m_q}+\sum_{i \neq q} b_{m_i}c_{m_{q-i}}=1+\sum_{i\neq q} b_{m_i}c_{m_{q-i}}
\]
hence all $c_{m_i}=0$ for $i\neq 0$. But this now means that there can be no double pairs, so all coefficients $b_i$ must be either $0$ or $1$.
\end{proof}
\begin{corollary}
If the factorisation of $f$ is unfair, then both $g$ and $h$ cannot be palindromic.
\end{corollary}
% ============================================================
\section{Large Gaps}
% ============================================================

We conclude with our main result on polynomials with large gaps. Motivated by work on the irreducibility of such polynomials~\cite{SawinShustermanStoll,FilasetaMatthews}, we consider the $0$--$1$ polynomial conjecture for polynomials of the form
\[
    f(x)=x^{n-\deg w}w(x)+v(x).
\]
The results of the previous sections allow us to obtain a sufficient condition under which such a polynomial is fair.

\begin{theorem}\label{thm:firstlargegap}
Let \(v(x)\) and \(w(x)\) be non-zero \(0\)--\(1\) polynomials and suppose $v(x)$ is fair. If
\[
    n\geq \max\{\deg w+2\deg v,\,2\deg w\}
\]
then the polynomial
\[
    f(x)=x^{\,n-\deg w}w(x)+v(x)
\]
is fair.
\end{theorem}

To prove the theorem, we begin with the following lemma.

\begin{lemma}\label{lm:aboveandbeyond}
Let \(v\), \(w\), and \(f\) be as in the statement of
Theorem~\ref{thm:firstlargegap}, and suppose
$ f=gh,$
where \(g\) and \(h\) have non-negative coefficients. Write
\(d_w=\deg w\). If
\[
    n\geq 2d_w,
\]
then either
\[
    b_i=0 \text{ for all } i\geq n-d_w
    \qquad\text{or}\qquad
    c_i=0 \text{ for all } i\geq n-d_w.
\]
\end{lemma}

\begin{proof}
Suppose that \(b_{n-j}>0\) for some \(j=0,1,\ldots,d_w\). We show
that
\[
    c_i=0 \qquad \text{for all } i\geq n-d_w.
\]

Suppose that \(c_{n-k}>0\) for some \(k=0,1,\ldots,d_w\). Then
\[
    b_{n-j}c_{n-k}>0,
\]
and hence
\[
    a_{2n-j-k}=1.
\]
Since \(n\geq 2d_w\),
\[
    2n-j-k\geq 2n-2d_w\geq n.
\]

If \(n>2d_w\), then \(2n-j-k>n\), which is impossible since
\(\deg f=n\).

It remains to consider the case \(n=2d_w\). Here
\(2n-j-k\leq n\) is possible only when \(j=k=d_w\). Thus the only
remaining possibility is
\[
    b_{d_w}>0
    \qquad\text{and}\qquad
    c_{d_w}>0.
\]
Since \(\deg g+\deg h=n=2d_w\), this forces
\(\deg g=\deg h=d_w\). Hence \(b_{d_w}\) and \(c_{d_w}\) are the
leading coefficients of \(g\) and \(h\), and therefore
\[
    b_{d_w}=c_{d_w}=1,
\]
contradicting Lemma~\ref{lm:no11}.

Thus, whenever some \(b_i\) with \(i\geq n-d_w\) is positive, all
\(c_i\) in this range are zero. If no such \(b_i\) is positive, then
\(b_i=0\) for all \(i\geq n-d_w\). Hence
\[
    b_i=0 \text{ for all } i\geq n-d_w
    \qquad\text{or}\qquad
    c_i=0 \text{ for all } i\geq n-d_w.
\]
\end{proof}

We now prove Theorem~\ref{thm:firstlargegap}.

\begin{proof}
Write \(d_v=\deg v\) and \(d_w=\deg w\).

Our goal is to show that there are no double pairs \((b_i,c_i)\) for
\(i=1,\ldots,n\). The result will then follow from Theorem~\ref{thm:doubletrouble}.

Since
\[
    n\geq \max\{d_w+2d_v,\,2d_w\},
\]
we have in particular \(n\geq2d_w\). By the preceding lemma,
\[
    b_i=0 \text{ for every } i\geq n-d_w
    \qquad\text{or}\qquad
    c_i=0 \text{ for every } i\geq n-d_w.
\]
In particular, there are no double pairs with index
\(i\geq n-d_w\).

We now consider the coefficients with indices \(i\leq d_v\). Since a
double pair at index \(i\) forces \(a_i=1\), no double pair can occur
at an index \(i\) for which \(a_i=0\). Thus it is enough to exclude
double pairs in the initial block.

By assumption, \(n\geq d_w+2d_v\). First suppose that
\[
    n>d_w+2d_v.
\]
Then \(n-d_w>2d_v\), so no term belonging to the shifted copy of
\(w(x)\) occurs in degrees \(0,\ldots,2d_v\). Hence
\[
    a_{d_v+i}=0, \qquad i=1,\ldots,d_v.
\]
Consequently, the coefficient equations
\[
    a_i=\sum_{j+k=i}b_jc_k,
    \qquad i=0,\ldots,2d_v,
\]
are precisely the coefficient equations for a factorisation of
\(v(x)\). Since \(v(x)\) is fair, no double pair can occur with
\(i\leq d_v\).

It remains to consider the boundary case
\[
    n=d_w+2d_v.
\]
Here the shifted copy of \(w(x)\) begins in degree \(2d_v\), so the
only equation in the reduced system which differs from that of
\(v(x)\) is the equation for \(a_{2d_v}\). Therefore the only
additional possibility which must be considered is that
\((b_{d_v},c_{d_v})\) is a double pair.

Suppose then that \(b_{d_v}>0\) and \(c_{d_v}>0\). Since
\(n-d_w=2d_v\), the preceding lemma implies that either
\[
    b_i=0 \text{ for every } i\geq2d_v,
\]
or
\[
    c_i=0 \text{ for every } i\geq2d_v.
\]
Moreover, since the coefficients between the two blocks of \(f\)
vanish, no positive coefficient of the corresponding factor can occur
between \(d_v\) and \(2d_v\). Thus either \(b_{d_v}\) or \(c_{d_v}\)
is the last positive coefficient of its respective factor. The leading
non-zero coefficient of each factor is equal to \(1\), whereas
\[
    0<b_{d_v}<1,
    \qquad
    0<c_{d_v}<1,
\]
because \((b_{d_v},c_{d_v})\) is a double pair. This is a
contradiction.

Hence no double pair can occur. By Theorem~\ref{thm:doubletrouble}, the
factorisation is fair.
\end{proof}

\begin{corollary}
If either $v(x)$ or $w(x)$ is fair and
\[
n\geq \deg(v)+\deg(w)+max\{deg(v),\deg(w)\}
\]
then $f$ is fair.
\end{corollary}

\begin{proof}
Since \(f\) is fair if and only if \(f^*\) is fair, and
\[
    f^*(x)
    =x^{\,n-\deg v}v^*(x)+w^*(x),
\]
it is enough to consider the case in which \(v(x)\) is fair; the
other case follows symmetrically from the fact that $w^*(x)$ is fair iff $w(x)$ is fair and our condition on $n$ is symmetric.

Now, notice that since $\deg(v)\geq 0$ we will have 
\[
\begin{aligned}
n
&\geq \deg(v)+\deg(w)+\max\{\deg(v),\deg(w)\} \\
&\geq \max\{\deg(w)+2\deg(v),\,2\deg(w)\}
\end{aligned}
\]
which is the exact condition required for Theorem~\ref{thm:firstlargegap} to hold.

\end{proof}

The result may be applied iteratively to construct families of fair 0--1 polynomials consisting of arbitrarily many sufficiently separated blocks.

\begin{question}
Can the bound on \(n\) be improved?
\end{question}

We know that the conjecture is known for palindromic
\(0\)--\(1\) polynomials. Therefore it is natural to consider
\[
    f(x)
    =x^{\,n-\deg v}v^*(x)+v(x).
\]
For sufficiently large \(n\), the low-degree reduced system associated
with \(f\) reproduces the corresponding system for \(v\).

This suggests trying to reverse the direction of the preceding theorem.
If such a converse held, then the palindromicity of \(f\) would imply
that \(f\) is fair, and one could hope to deduce that \(v\) itself is fair.
Since \(v\) is arbitrary, such a result would resolve the conjecture.

\begin{question}
Can we prove a converse? More precisely, if
\[
    f(x)=x^{\,n-\deg w}w(x)+v(x)
\]
is fair, under what additional hypotheses must \(v\) and \(w\) also be
fair?
\end{question}

\paragraph{Acknowledgements.}
I would like to thank Professor Samir Siksek for his supervision, helpful
discussions, and feedback throughout this project.

\end{document}